\documentclass[11pt]{amsart}
\usepackage{amsfonts,amssymb,amscd}
\usepackage[all]{xy}

\newtheorem{lemm}{Lemma}[section]
\newtheorem{theo}[lemm]{Theorem}

\newtheorem{prop}[lemm]{Proposition}
\newtheorem{defi}[lemm]{Definition}

\def\SS{{\mathbb S}}
\def\OP2{\mathbb{OP}^2}
\def\CC{{\mathbb C}}
\def\PP{{\mathbb P}}
\def\QQ{{\mathbb Q}}
\def\ZZ{{\mathbb Z}}
\def\ra{\rightarrow}
\def\cO{\mathcal{O}}\def\cD{\mathcal{D}}\def\cA{\mathcal{A}}
\def\cJ{\mathcal{J}}
\def\cX{\mathcal{X}}
\def\cB{\mathcal{B}}
\def\fP{\mathfrak{P}}

\def\cX{\mathcal{X}}

\def\cD{\mathcal{D}}
\def\cE{\mathcal{E}}
\def\lra{\longrightarrow}

\begin{document}

\title{Fano manifolds of Calabi-Yau type}
\author[A. Iliev]{Atanas Iliev}
\address{Department of Mathematics,
Seoul National University,
Seoul 151-747, Korea}
\email{{\tt ailiev2001@yahoo.com}}
\author[L. Manivel]{Laurent Manivel}
\address{Institut Fourier,  
Universit\'e de Grenoble et CNRS,
BP 74, 38402 Saint-Martin d'H\`eres, France}
\email{{\tt Laurent.Manivel@ujf-grenoble.fr}}

\begin{abstract}
We introduce and we study a class of odd dimensional compact 
complex manifolds whose Hodge structure in middle dimension
looks like that of a Calabi-Yau threefold. We construct 
several series of interesting examples from rational homogeneous 
spaces with special properties. 
\end{abstract}

\maketitle

\section{Introduction}

If we think, quite naively, to mirror symmetry as a kind of mysterious 
process exchanging families of Calabi-Yau threefolds in such a way that 
the Hodge diamond rotates by a quarter-turn, an obvious difficulty 
arises for {\it rigid} Calabi-Yau's: their mirrors could not be K\"ahler!
The solution given to that issue in \cite{cdp, sch} 
was to describe the mirror of certain rigid
Calabi-Yau threefolds as families of higher dimensional manifolds of
a very special kind. In particular, these mirror manifolds have odd dimension
$2n+1$ and their Hodge structure in middle dimension looks like that
of a Calabi-Yau threefold, in the sense that the only non-zero Hodge 
numbers are $h^{n-1,n+2}=h^{n+2,n-1}=1$ and $h^{n,n+1}=h^{n+1,n}$. 

This is precisely the kind of manifolds that we study in this paper, 
with the idea that, under some conditions that will be made more precise 
below, these {\it manifolds of Calabi-Yau type} should share some of the 
very nice properties of Calabi-Yau threefolds. The definition that we 
will use  may not exactly be the correct one and we consider it as 
provisional. Our main general result will be that for a  
manifold of Calabi-Yau type with non-obstructed deformations,
it implies  that the relative 
intermediate Jacobian forms an integrable system over the gauged 
moduli space -- a result due to Donagi and Markman for Calabi-Yau 
threefolds. The period map also has a very similar behavior to the 
Calabi-Yau setting. Related homological properties are explored in \cite{ik},
where it is proved that one can extend at least to certain Fano manifolds 
of Calabi-Yau type, the celebrated result of Voisin according to which 
the Griffiths group of a Calabi-Yau threefold cannot be finitely generated.

But the main goal of the present paper is to construct examples 
of manifolds of Calabi-Yau type, which will all be Fano's. 
Our main source of examples will be complete intersections in 
homogeneous spaces: we will show that, for such a  construction 
to work, one needs strong numerical coincidences, notably between
the dimension and the index of the ambient homogeneous space. 
These coincidences are observed in some cases. Notably, we construct
interesting examples as quadratic sections of homogeneous spaces that 
are Mukai varieties of even dimension (Proposition \ref{muk-pair}),
while for Mukai varieties of odd dimension we need to take double covers
branched over quadratic sections  (Proposition \ref{muk-impair}). 

But our most intriguing series of examples is constructed from
homogeneous spaces with the property that their projective dual
is a hypersurface of degree equal to the coindex minus one (the 
coindex of a Fano manifold being defined as the dimension minus
the index). We are aware of four cases for which this strange
coincidence can be observed, and we show that suitable linear 
sections provide examples of Fano manifolds with the required
Hodge numbers. In fact we can really conclude that they are 
of Calabi-Yau type only in half of the cases; for the two others
there remains a tedious computation to be done, but we have no
doubt about the fact that the final conclusion should be positive. 
In each case, we would conclude that a certain type of hypersurfaces
can be, generically,  represented as a linear section of the dual
hypersurface to our homogeneous space, in a finite number of ways. 
For example, a generic cubic sevenfold can be represented as 
a linear section of the famous Cartan cubic, the $E_6$-invariant
hypersurface in $\PP^{26}$, in a finite number of different ways.
This has interesting consequences for its derived category, 
as we show in \cite{im-cub} where this example in studied in more 
details. 

It would be interesting to find out more examples of Fano 
manifolds of Calabi-Yau type. A possible strategy would be to
start from rigid Calabi-Yau threefolds and try to identify 
their mirrors systematically. We expect that these varieties 
will exhibit a rich and interesting geometry.

\section{Definition and first properties}

We start with our main definition. 

\begin{defi}
Let $X$ be a smooth complex compact variety of  odd dimension $2n+1$, $n\ge 1$. 
We call $X$ a {\it manifold of Calabi-Yau type} if
\begin{enumerate}
\item
The middle dimensional Hodge structure is numerically similar to that of a 
Calabi-Yau threefold, that is 
$$h^{n+2,n-1}(X)=1, \qquad and\quad h^{n+p+1,n-p}(X)=0 \;\;for\;p\ge 2.$$
\item
For any generator $\omega\in H^{n+2,n-1}(X)$, the contraction map 
$$H^1(X,TX)\stackrel{\omega}{\lra} H^{n-1}(X,\Omega_X^{n+1})$$
is an isomorphism. 
\item
The Hodge numbers $h^{k,0}(X)=0$ for $1\le k\le 2n$.
\end{enumerate}
\end{defi}

A Calabi-Yau threefold is of course a manifold of Calabi-Yau type. 
The definition may not be the optimal one, since in particular, 
even in dimension three, there exist manifolds of Calabi-Yau type
which are not Calabi-Yau stricto sensu. Also, the condition of 
being a manifold of Calabi-Yau type is clearly invariant under
small deformations but probably not under arbitrary deformations, 
contrary to being Calabi-Yau. 

Rather, our main motivation is to find interesting examples of 
manifolds of Calabi-Yau type, and to investigate their geometry. 
In many respects they will behave like  Calabi-Yau threefolds,
and sometimes the examples we will meet will be related, in a 
rather non trivial way, to Calabi-Yau threefolds. 
Sometimes we will observe interesting differences of behavior. 

Deformations of Calabi-Yau manifolds are always unobstructed, 
as follows from the celebrated Tian-Todorov theorem. We don't know 
whether this always remains true for manifolds of Calabi-Yau type.
(Actually all the concrete examples that we deal with are Fano,
hence have  unobstructed deformations.)
Nevertheless, our first observation is that versal deformations 
of manifolds of Calabi-Yau type give rise, exactly as families 
of Calabi-Yau threefolds do, to some beautiful integrable systems. 
 
Consider a versal family $\pi :\cX\ra\cB$ of manifolds of Calabi-Yau
type, with special fiber $X$. We suppose that the base $\cB$ is
 identified with  an open subset of  $H^1(X,TX)$. Since $h^{n+p+1,n-p}(X)=0$ for
$p\ge 2$, the line bundle $R^{n-1}\pi_*\Omega_{\cX/\cB}^{n+2}$ over $\cB$ 
is holomorphic, and the complement $\cB_*$ of the zero section, with the 
pull-back $\cX_*$ of $\cX$, defines the family of {\it gauged} manifolds
of  Calabi-Yau type. 

Recall that the intermediate Jacobian of $X$ is the complex torus 
$$J(X)=F^{n+1}H^{2n+1}(X,\CC)^\vee/ H_{2n+1}(X,\ZZ).$$
Since $X$ is of Calabi-Yau type, the dimension of $J(X)$ is 
$h^{n+1,n}(X)+1=h^1(TX)+1=\dim \cB_*$. Globalizing the construction,
we get a torus bundle, the relative  intermediate Jacobian
$$\cJ(\cX/\cB)\lra  \cB_*.$$
The theorem proved by Donagi and Markman  for Calabi-Yau threefolds 
\cite{dm}
can be extended to our generalized setting:

\begin{theo}\label{int} 
The torus fibration $\cJ(\cX/\cB)\lra  \cB_*$ defines a completely integrable
Hamiltonian system. 
\end{theo}

\proof The proof of the Theorem for Calabi-Yau threefolds is easy to extend to our 
setting. The main observation is that the tangent exact sequence of the 
$\CC^*$-bundle $\cB_*\ra \cB$ can be identified with the Hodge filtration 
of $H^{2n+1}(\cX_*/\cB_*,\CC)$. More precisely  there is a commutative diagram
$$\begin{array}{ccccccc}
0\lra &  T_{(X,\omega)}(\cX_*/\cX) & \lra & T_{(X,\omega)}\cX_* &
\lra & T_X\cX & \lra 0 \\
 & \downarrow & & \downarrow & & \downarrow & \\
 0\lra &  H^{n+2,n-1}(X) & \lra & F^{n+1}H^{2n+1}(X,\CC) &
\lra & H^{n+1,n}(X) & \lra 0,
\end{array}$$
where the vertical maps are all isomorphisms. In particular  the middle one is 
given by the differential of $\omega$ (considered as a section on $\cX_*$ of the 
holomorphic  bundle $F^{n+2}H^{2n+1}(\cX_*/\cB_*,\CC)$), by the Gauss-Manin
connection: this yields, by Griffiths' transversality, a map
$$\nabla\omega : T_{(X,\omega)}\cX_*\lra F^{n+1}H^{2n+1}(X,\CC),$$
which must be an isomorphism since its two graded parts are. 

This observation yields a global identification between the cotangent bundle 
of $\cX_*$ and the dual of the Hodge bundle $F^{n+1}H^{2n+1}(\cX_*/\cB_*,\CC)$. 
Now $\Omega_{\cX_*}$ has a canonical symplectic structure. We need to check that it 
descends to the torus fibration $\cJ(\cX/\cB)$, which is just a quotient by the 
relative lattice defined by locally constant $(n+1)$-cycles. So there  remains 
to verify that such a locally constant cycle, say $\gamma$, when considered 
as a one-form on $\cX_*$, is closed. But the argument of Donagi-Markman, which 
consists in constructing a local primitive by integrating $\omega$ over $\gamma$,
adapts mutatis mutandis.
\qed

\smallskip\noindent {\it Remark}. This implies the existence of a special K\"ahler 
structure on the base, with an extremely rich geometry \cite{fr}. 
 
\medskip
The existence of that integrable system is related to that of a certain 
cubic form on the base manifold, which in the usual Calabi-Yau setting is 
the famous Yukawa cubic. The Yukawa cubic is defined, for a gauged Calabi-Yau threefold
 $(X,\omega)$, by the composition 
$$Sym^3H^1(X,TX)\lra H^3(X,\omega_X^{-1})\stackrel{SD}{\lra} H^0(X,\omega_X^2)^\vee
\stackrel{\omega^2}{\lra}\CC,$$
where $SD$ denotes Serre duality. In our generalized setting, this construction 
can be extended as follows. First consider the natural map $\wedge^{n-1}TX\otimes 
\wedge^{n-1}TX\ra \wedge^{2n-2}TX$, which is symmetric if $n$ is odd and skew-symmetric
if $n$ is even. Twisting by $\omega_X^2$ we get a map $\Omega_X^{n+2}\otimes\Omega_X^{n+2}
\ra\Omega_X^3\otimes\omega_X$ with the same parity. We get an induced map (always
symmetric, since the cup-product has the same parity as the degree, here $n-1$) 
$$Sym^2H^{n-1}(X,\Omega_X^{n+2})\lra H^{2n-2}(X,\Omega_X^3\otimes\omega_X),$$
sending the square of a generator $\omega$ of $H^{n-1}(X,\Omega_X^{n+2})$ to 
an element $s(\omega)$ of $H^{2n-2}(X,\Omega_X^3\otimes\omega_X)$. 
We can then define a map
$$Sym^3H^1(X,TX)\lra H^3(X,\wedge^3TX)\stackrel{s(\omega)}{\lra} 
H^{2n+1}(X,\omega_X)=\CC,$$
which is our generalized cubic at the point of $\cB_*$ defined by the gauged manifold
of Calabi-Yau type $(X,\omega)$. 

\medskip 
For a slightly different view-point, consider the period map 
$$\begin{array}{rccl}
\fP : & \cB_* &\ra & H^{2n+1}(X_0,\CC),\\
 & (X,\omega) & \mapsto &\omega,
\end{array}$$
defined once the family $\cX\ra\cB$ has been trivialized (in the differentiable
category). The same proof as in the usual Calabi-Yau setting leads to the following 
statement:

\begin{prop}
The period map of a versal family of manifolds of Calabi-Yau type is a local immersion.
  
Moreover, its image is a pointed cone in $H^{2n+1}(X_0,\CC)$,
which is Lagrangian relatively to the symplectic structure defined by
the intersection product.  
\end{prop}

It would be interesting to investigate under which conditions a Frobenius
structure can be defined from our generalized Yukawa cubic. 
We know from \cite{cdp} that everything works fine for one of the most 
interesting examples of manifolds of Calabi-Yau type, the cubic sevenfold. 
Is there any chance to extend this to other examples? One can make a 
first step in this direction by mimicking the proof of Proposition 3.3 
in \cite{voisin}, which allows one to  prove that:

\begin{prop}
One can choose a local section $\omega$ of $F^{n+2}H^{2n+1}(\cX_*/\cB_*,\CC)$
 in such a way that the Yukawa cubic 
$$Y^{\omega} : Sym^3H^1(X,TX)\lra  H^3(X,\wedge^3TX)\stackrel{s(\omega)}{\lra}\CC$$
is defined by the derivatives of some potential $F$ on $\cB$.  
\end{prop}

Understanding under which conditions this potential would satisfy the WDVV equation
seems to be a very difficult problem, but the question is posed: could there
be a version of mirror symmetry for (some) manifolds of Calabi-Yau type?

\section{Complete intersections of Calabi-Yau type}

 \subsection{Hodge numbers of complete intersections}

In this section we review well-known results about the cohomology 
of complete intersections in weighted projective spaces. Let us begin 
with a degree $d$ hypersurface $X_d\subset w\PP^n$, where 
$w=(w_0,\ldots ,w_n)$ is the set of weights defining the weighted 
projective space $w\PP^n$. 

The usual projective space $\PP^n$ corresponds to $w=(1,\ldots ,1)$.
In that case, it is well known that the primitive cohomology of $X_d$
can be expressed in terms of the Jacobian ring
$$R=\CC[x_0,\ldots ,x_n]/(\frac{\partial F}{\partial x_0},\ldots 
\frac{\partial F}{\partial x_n}),$$
where $F$ denotes an equation of $X_d$. In fact there exists a general
relation of that kind for any quasi-smooth  hypersurface 
$X_d\subset w\PP^n$, which is due to Steenbrink \cite{steenbrink}:
$$H^{n-p-1,p}(X_d)_0\simeq R_{(p+1)d-|w|},$$
where $|w|=w_0+\cdots +w_n$ and $R_k$ denotes the degree $k$ component
of $R$, with respect to the natural grading defined by $\deg (x_i)=w_i$. 

An obvious consequence is that if $d$ divides $|w|$, that is, $|w|=
(p+1)d$ form some $p$, then 
\begin{eqnarray*}
 &H^{n-q-1,q}(X_d)_0=0\quad \mathrm{for}\; q\le p-1, \\
 &H^{n-p-1,p}(X_d)_0\simeq R_0\simeq \CC, \\
 &H^{n-p-2,p+1}(X_d)_0\simeq R_{d}\simeq H^1(X_d,TX_d).
\end{eqnarray*}

In particular, if $n=2m$ is even and $|w|=(m-1)d$, the primitive
cohomology of $X_d$ defines a Hodge structure of the same type as that
of a Calabi-Yau threefold. We find a short list of smooth examples 
in dimension $n-1\ge 4$:

$$\begin{array}{ccccl}
n-1 & w & d & moduli & X\\
7 & 1^9 & 3 & 84 & \mathrm{cubic\;sevenfold} \\
5 & 1^7 & 2,3 & 83 & \mathrm{cubic\;section\;of\;\QQ^6} \\
5 & 1^6,2 & 4 & 90 & \mathrm{double\;quartic\;fivefold}
\end{array}$$

\smallskip
We will meet again the cubic sevenfolds later on. They 
will be investigated in a more systematic way in \cite{im-cub}. 

\subsection{A technical lemma}

We would like to generalize the previous examples and show that certain 
low degree hypersurfaces in homogeneous spaces are of Calabi-Yau type. 

The  general set-up would be the following. Suppose that $\Sigma$ is a Fano variety 
of dimension $d_\Sigma$, with $-K_\Sigma = \iota L$ for some ample line bundle $L$. 
Suppose moreover that $L$ is generated by global sections and consider the zero-locus
$X$ of a general section of $L$. Its dimension is of course $d_X=d_\Sigma-1$, and 
$-K_X = (\iota-1) L_{|X}$. The proof of our next technical lemma is a rather standard
play with long exact sequences. It is essentially an application of Griffiths' techniques,
for the computation of Hodge structures of hypersurfaces \cite{Gr, voisin}.

\begin{lemm}\label{tech}
Let $\Sigma$ be such that $d_\Sigma =2\iota+2$. 
\begin{enumerate}
\item Suppose that for $0\le k<p<\iota$,
$$H^{d_\Sigma-p+k-1}(\Sigma,\Omega_{\Sigma}^{p-k}\otimes L^{-k})=
H^{d_\Sigma-p+k}(\Sigma,\Omega_{\Sigma}^{p-k}\otimes L^{-k-1})=0.$$
Then $h^{p,d_X-p}(X)=0$ for $p<\iota$.
\item Suppose that for $1\le k<\iota$ and $\epsilon,\epsilon'\in\{0,1\}$, 
$$H^{\iota+k+1+\epsilon+\epsilon'}(\Sigma,\Omega_{\Sigma}^{i-k}\otimes L^{-k+1-\epsilon})=0.$$
Then $h^{\iota+2,\iota-1}(X)=1$.
\item Suppose moreover that for $2\le k\le\iota$ and $\epsilon,\epsilon'\in\{0,1\}$, 
$$H^{\iota+k+1+\epsilon+\epsilon'}(\Sigma,\Omega_{\Sigma}^{i-k}\otimes L^{-k-\epsilon})=0.$$
\end{enumerate}
Then $X$ is a manifold of Calabi-Yau type. 
\end{lemm}

\proof Use the conormal sequence of the pair $(X,\Sigma)$. Its $p$-th wedge power is 
the complex 
$$0\ra L^{-p}_{|X}\ra \Omega_\Sigma^1\otimes L^{-p+1}_{|X}\ra\cdots 
\ra\Omega_\Sigma^{p-1}\otimes L_{|X}\ra\Omega^p_{\Sigma|X}\ra\Omega_X^p\ra 0.$$
Under hypothesis (1), for $p\le\iota-2$, we have 
$$H^{d_X-p+k}(X,\Omega_\Sigma^{p-k}\otimes L^{-k}_{|X})=0$$
for all $k$ from $0$ to $p$. Our first claim follows. 

For $p=\iota-1$, our hypothesis imply that 
$$H^{d_X-p+k+\epsilon}(X,\Omega_\Sigma^{p-k}\otimes L^{-k}_{|X})=0$$
for $\epsilon\in\{0,1\}$. But then, it follows that 
$$H^{\iota+2}(X,\Omega_X^{\iota-1})=H^{d_X}(X,L^{-\iota+1}_{|X})=H^{d_X}(X,K_X)=\CC,$$
which was our second claim. 

Finally, for $p=\iota$ the hypothesis we made are such that we get a complex
$$0\ra H^{\iota+1}(X,\Omega_X^{\iota})\ra  H^{d_X}(X,L^{-\iota}_{|X})
\ra  H^{d_X}(X,\Omega_\Sigma^1\otimes L^{-\iota+1}_{|X}).$$
Using Serre duality, we get the dual complex
$$H^0(X,T\Sigma_{|X})\ra H^0(X,L)\ra H^{\iota}(X,\Omega_X^{\iota+1})\ra 0.$$
But the cokernel of the first arrow is precisely $H^1(X,TX)$, at least if 
we have $H^1(X,T\Sigma_{|X})=0$ or, by Serre duality again, 
$H^{d_X-1}(X,\Omega_\Sigma^{d_X}\otimes L^{-\iota}_{|X})=0$. 
This follows from (3). \qed 

\medskip
Of course we can easily imagine variants if this lemma for complete intersections, or even 
more generally, for zero-loci of global sections of vector bundles. For the construction 
to work, we need a huge number of vanishing conditions on the ambient variety $\Sigma$, 
and it is not so clear a priori that we can find any example that way. We have been
able to find some under the hypothesis that $\Sigma$ is a homogeneous space.

\subsection{Linear sections of  homogeneous spaces}

Suppose that $\Sigma$ is a rational homogeneous variety, and that $\iota=\iota_\Sigma$
is equal to its index. In particular we require that the dimension of $\Sigma$ is $d_\Sigma
=2\iota_\Sigma+2$. There is a short list of suitable examples:

 $$\begin{array}{cccc}
\Sigma & dim & index & moduli\\
(\PP^1)^6 & 6 & 2 & 45 \\
(\PP^1)^3\times\PP^3 & 6 & 2 & 55 \\
(\PP^2)^4 & 8 & 3 & 48 \\
 (\PP^4)^3 & 12 & 5 & 52 \\
G(2,5)\times G(2,5) & 12 & 5 & 51 \\
G(4,9) & 20 & 9 & 46\\
G(3,11) & 24 & 11 & 45
\end{array}$$

\medskip

\begin{prop}
Let $\Sigma$ belong to the list above, and $X=\Sigma\cap H$ be a smooth hyperplane section of 
$\Sigma$ in its minimal complete embedding. 
Then $X$ is of Calabi-Yau type. 
\end{prop}

\proof We check case by case that the vanishing conditions of Lemma \ref{tech} do hold.  
For Grassmannians this requires the diagrammatic methods first developed in \cite{snow}.
\qed 

\subsection{Quadratic sections of  homogeneous spaces}

Now suppose that $\Sigma$ is again a  rational homogeneous variety of even index
$\iota_\Sigma$, and let $\iota=\iota_\Sigma/2$. The relation $d_\Sigma=2\iota+2$
becomes $\iota_\Sigma=d_\Sigma-2$, which is the definition of the {\it Mukai varieties}. 
It has been proved that apart from complete intersections, Mukai 
varieties are linear sections of homogeneous varieties (note that
a linear section of a Mukai variety, if Fano, is again a Mukai variety), 
or quadric sections of the cone over $G(2,5)$. 
The homogeneous Mukai varieties are the following:
 $$\begin{array}{cccc}
\Sigma & dim & index & moduli\\
(\PP^1)^4 & 4 & 2 & 68 \\ 
G_2^{ad} & 5 & 3 & 62 \\
\PP^3\times \PP^3 & 6 & 4 & 69\\
LG(3,6) & 6 & 4 & 62\\
IG(2,6) & 7 & 5 & 68\\
G(2,6) & 8 & 6 & 69\\
\SS_{10} & 10 & 8 & 80 
\end{array}$$

In this table we have denoted by $IG(k,2n)$ the symplectic  Grassmannian,
parametrizing isotropic $k$-dimensional subspaces of a $2n$-dimensional vector space
endowed with a non-degenerate skew-symmetric form. When $k=n$ we rather 
use the notation $LG(n,2n)$ and call this variety the Lagrangian Grassmannian. 
Similarly, $\SS_{2n}$ denotes the spinor variety, which parametrizes one of the
two families of isotropic $n$-dimensional subspaces of a $2n$-dimensional 
vector space endowed with a non-degenerate quadratic form. The terminology
comes from the fact that the minimal equivariant embedding of this variety is
inside a projectivized half-spin representation -- a quadratic Veronese is 
required in order to recover the Pl\"ucker embedding. Finally, the 
adjoint variety $G_2^{ad}$ denotes the closed $G_2$-orbit inside the 
projectivized adjoint representation.

\begin{prop}\label{muk-pair}
Let $X=\Sigma\cap Q$ be a smooth quadric section of a homogeneous Mukai 
variety of even dimension. Then $X$ is of Calabi-Yau type. 
\end{prop}

 \proof 
Again we need to check, case by case, the conditions of Lemma \ref{tech}.
Let us treat in some detail the case of $\SS_{10}$, which is the most 
complicated one. 
Here $X$ is nine dimensional, and the claim is that its middle dimensional 
Hodge numbers are 
$$h^{9,0}=h^{8,1}=h^{7,2}=0, \quad h^{6,3}=1, \quad h^{5,4}=80.$$
Let $E$ denote the tautological rank five vector bundle on $\Sigma=\SS_{10}$. 
The cotangent bundle of $\Sigma$ is $\Omega_{\Sigma}^1\simeq\wedge^2E$. 
In higher degree, the bundle $\Omega_{\Sigma}^p$ of $p$-forms can have
several components, which are Schur powers of $E$ obtained as follows. 
Let $a=(a_0>a_1>\cdots >a_m)$ be a decreasing sequence of integers, whose sum
is equal to $p$. Define the sequence 
$$\lambda(a)=(a_0,a_1+1,\ldots,a_m+m,m^{a_{m-1}-a_m-1},
\ldots ,1^{a_0-a_1-1}).$$
Then by the results of \cite{kos}, 
$\Omega_{\Sigma}^p=\wedge^p(\wedge^2E)$ is the direct sum of the 
Schur powers $S_{\lambda(a)}E$. Note that since $E$ has rank five, 
the number of non zero components of $\lambda(a)$ must not exceed 
five, which means that $a_0-a_m\le 4$.

Bott's theorem allows to compute the cohomology groups of a Schur power
$S_{\lambda}E$ of the tautological bundle, for $\lambda=(\lambda_1,\ldots,
\lambda_5)$ a non increasing sequence of non negative numbers. The rule 
is the following. Suppose that there exists a pair $(i,j)$ of distinct indices
such that $\lambda_i+\lambda_j=i+j-2$. Then $S_{\lambda}E$ is acyclic. 
Otherwise, let $q$ be the number of pairs such that $\lambda_i+\lambda_j 
>i+j-2$. Then $S_{\lambda}E$ has a unique non-zero cohomology group,
that of degree $q$. 

Applying these rules, it is easy to check the following 
statement. 

\begin{lemm}
The following bundles of twisted forms on $\Sigma$ are acyclic:
\begin{enumerate}
\item $\Omega_{\Sigma}^1(-k)$ for $k=1,2,3$.
\item $\Omega_{\Sigma}^2(-k)$ for $k=1,2,3$.
\item $\Omega_{\Sigma}^3(-k)$ for $k=1,2$.
\item $\Omega_{\Sigma}^4(-k)$ for $k=1,2$.
\end{enumerate}
Moreover, for $k\ge 4$, $\Omega_{\Sigma}^1(-k)$ and $\Omega_{\Sigma}^2(-k)$
have non-zero cohomology groups only in degree $\dim\Sigma$. And the same
is true for $\Omega_{\Sigma}^3(-k)$ and $\Omega_{\Sigma}^4(-k)$ when 
$k\ge 3$.
\end{lemm}

This is precisely what we need in order to apply Lemma \ref{tech}. \qed

\subsection{A quasi-homogeneous Mukai variety}

Consider a hyperplane section $\Theta$ of the Lagrangian Grassmannian
$LG(3,6)$. This is a five dimensional variety of index three.

\begin{prop} The variety $\Theta$ has the following properties:
\begin{enumerate}
\item it is a quasi-homogeneous $PSL_3$-variety;
\item its non-zero Hodge numbers are $h^{p,p}(\Theta)=1$, $1\le p\le 5$;
\item it is a minimal compactification of $\CC^5$;
\item it is rigid.
\end{enumerate}
\end{prop}

\proof Consider the Lagrangian Grassmannian $LG(3,6)$ in its minimal 
homogeneous embedding $\PP V=\PP^{13}$. The orbit structure of $Sp_6$ in 
$\PP V$ is well-known, in particular it is prehomogeneous: there exists an
open orbit, and the generic stabilizer is $PSL_3$ (up to a finite group). 
Since $V$ is self-dual, we conclude that the general hyperplane section 
$\Theta$ is rigid and admits a $PSL_3$-action. 

An explicit computation shows that $PSL_3$ has four orbits in $\Theta$. 
This can be seen as follows. The Lagrangian Grassmannian $LG(3,6)$ is known
to be one of the few homogeneous varieties with one apparent double point. 
This means that a general point of $\PP V$ belongs to a unique secant 
$\overline{xy}$, where $x$ and $y$ represent two generic isotropic planes. 
The orbits in $\Theta$ are then determined by the relative position 
with respect to $x$ and $y$, that is, by the dimensions of these intersections.
These can be $(0,0)$ (the generic case, giving the open orbit), $(1,1)$ 
(a hyperplane section), $(2,1)$ or $(1,2)$ (two Veronese surfaces). 

The computation of the Hodge numbers is straightforward. Finally, it 
is easy to check that a maximal torus in $PSL_3$ has a finite number
of fixed points in $\Theta$. The Byalinicki-Birula decomposition \cite{bb}
then ensures that $\Theta$ is a compactification of $\CC^5$, and it is clearly 
minimal. \qed

\subsection{Double coverings}

Consider a double covering $Y\ra \Sigma$, branched over a smooth hypersurface
$X$ of $\Sigma$, of degree $2d$. The Hodge numbers of $Y$ can be computed 
in terms of the pair $(\Sigma,X)$ as follows \cite{cynk}:
 $$h^{p,q}(Y)=h^{p,q}(\Sigma)+h^q(\Omega^p_\Sigma(\log X)(-d)).$$
Recall that $\Omega^p_\Sigma(\log X)$ denotes the vector bundle of 
logarithmic $p$-forms on $\Sigma$ with simple poles at most on $X$. 
There is an exact sequence 
$$0\ra \Omega^p_\Sigma\ra \Omega^p_\Sigma(\log X)\ra \Omega^{p-1}_X\ra 0$$
defined by taking residues on $X$. 

Suppose that $\Sigma$ is homogeneous of odd dimension $n=2m+1$. 
In particular the Hodge structure of $\Sigma$ is pure and we get 
$h^{p,q}(\Sigma)=0$ for $p+q=n$. Suppose moreover that 
$H^q(\Sigma,\Omega^p_\Sigma(-d))=
H^{q+1}(\Sigma,\Omega^p_\Sigma(-d))=0$. Then $h^{p,q}(Y)=
h^q(X,\Omega^{p-1}_X(-d))$, 
which can be computed by considering a long exact sequence like
in the proof of Lemma \ref{tech}. Under favorable circumstances, we expect to 
derive that $$h^{m+2,m-1}(Y)=
h^{m+2}(X,\Omega^{m-2}_X(-d))=h^{2m}(X,\cO_X(-d-(m-2)2d))=\CC$$
if moreover $\omega_X=\cO_X(-(2m-3)d)$, that is, $\omega_\Sigma=
\cO_\Sigma(-(2m-1)d)$. If we let $d=1$, this means that $\Sigma$
is a Mukai variety. If we let $\Sigma$ be homogeneous, or the 
quasi-homogeneous Mukai variety of the previous section, 
we check that the expected vanishing theorems do hold and we arrive
at the following conclusion. 

\begin{prop}\label{muk-impair}
Let $Y\ra\Sigma$ be a double covering of a homogeneous or quasi-homogeneous 
Mukai variety of odd dimension, branched over a smooth quadric section. 
Then $Y$ is a Fano manifold of Calabi-Yau type. 
\end{prop}

\subsection{The Cayley trick}
Let $X\subset\PP^N$ be a complete intersection 
of multidegree $(d_1,\ldots ,d_m)$, where none of these degrees equals one. 
Let $L_i=\cO_{\PP^N}(d_i)$ and $P=\PP(E^\vee)$,
where $E=L_1\oplus\cdots\oplus L_m$. The equations of $X$ can be considered as a section 
$\sigma$ in $H^0(\PP^N,E)\simeq H^0(P,\cO_E(1))$. Hence an associated divisor $\cX\subset P$,
which is smooth when $X$ is smooth. By \cite[Lemma 2.7]{nagel}, the cohomology of $\cX$ 
is essentially the same as that of $X$, with a shift by $2m-2=\dim\cX-\dim X$ in the degrees.

Moreover the deformations of $X$ and $\cX$ are both unobstructed, and
\begin{equation}\label{cayley}
H^1(X,TX)\simeq H^1(\cX,T\cX).
\end{equation}
Indeed, the left hand side of (\ref{cayley}) is given by the exact sequence
$$ H^0(X,T\PP^N_{|X})\ra  H^0(X,E_{|X})\ra H^1(X,TX)\ra 0.$$
It is easy to check that $ H^0(X,T\PP^N_{|X})\simeq  H^0(\PP^N,T\PP^N)=\mathfrak{sl}_N$. 
Moreover, the Koszul complex of the section $\sigma$ of $E$ defining $X$ shows that $H^0(X,E_{|X})$
is just the cokernel of the natural map
$$H^0(\PP^N,End(E))\ra H^0(\PP^N,E)$$
defined by applying an endomorphism of $E$ to the section $\sigma$.  

Similarly, the right hand side of (\ref{cayley}) is given by an exact sequence
$$ H^0(\cX,TP_{|\cX})\ra  H^0(\cX,\cO_{\cX}(1))\ra H^1(\cX,T\cX)\ra 0.$$
We have $H^0(\cX,\cO_{\cX}(1))\simeq H^0(P,\cO_E(1))/\CC\sigma=H^0(\PP^N,E)/\CC\sigma$. 
Moreover, $H^0(\cX,TP_{|\cX})=H^0(P,TP)$ can be computed through the differential 
of the projection $\pi$ from $P$ to $X$ and the description of its kernel, 
the vertical tangent space $T^vP$, by the relative Euler sequence
$$0\ra \cO_P\ra  \pi^*E^\vee\otimes\cO_P(1)\ra T^vP\ra 0.$$
We get that $H^0(P,TP)$ is the direct sum of $H^0(P,T^vP)$ and $\mathfrak{sl}_N$, while 
$H^0(P,T^vP)\simeq H^0(\PP^N,End_0(E))$. Since the maps $H^0(\PP^N,End(E))\ra H^0(\PP^N,E)$
and  $H^0(\PP^N,End_0(E))\ra H^0(\PP^N,E)/\CC\sigma$ have the same cokernels, 
the identification (\ref{cayley}) follows. We conclude:

\begin{prop} If $X\subset\PP^N$ is a Calabi-Yau manifold, or a complete intersection of
Calabi-Yau type, then $\cX\subset P$ is also of Calabi-Yau type.
\end{prop}

\section{A special series of manifolds of Calabi-Yau type}

\subsection{A special non-vanishing}
 
Consider the following series of homogeneous varieties $\Sigma$:

$$\begin{array}{ccccc}
\Sigma  & dim & index & coindex & deg \\
\OP2 & 16 &  12 & 4 & 3  \\
\SS_{12} & 15 &  10 & 5 & 4  \\
G(2,10) & 16 & 10 & 6 & 5 \\
\SS_{14} & 21 & 12 & 9 & 8 
\end{array}$$

\medskip
Here we denoted by $\OP2$ the so-called {\it Cayley plane}, 
which is the unique Hermitian symmetric space of type $E_6$. Recall
that the minimal representation of  $E_6$ has dimension $27$. 
They Cayley plane can be obtained as the minimal $E_6$-orbit
in the projectivization of this representation (or its dual). 
For more on the extraordinary properties of this variety, 
see \cite{im-cayley} and the references therein.

\medskip
Recall that the index $r$ is defined by the identity $K_{\Sigma}=-rH$, 
where $H$ is the ample generator of the Picard group of $\Sigma$. 
The coindex is the difference between the dimension and the index. 

The last column gives the degree of the dual variety of $\Sigma$. 
Note that this degree is always one less than the coindex, a 
coincidence that will play a major role below.  

\begin{lemm}\label{nonvan} 
Let $\Sigma$ be as above, with index $r$ and coindex $c$. 
Then 
$$H^{r+2}(\Sigma,\Omega_\Sigma^{c-2}(c-r-1))=\CC.$$
\end{lemm}

\proof This is a straightforward application of the Borel-Weil-Bott
theorem. The details are provided in the Appendix. \qed

\subsection{Relation with projective duality}

The explanation for the non-vanishing of Lemma \ref{nonvan} is the following. 
We start by recalling the usual setting of projective duality, 
which we summarize in the diagram

\centerline{
\xymatrix{ & & I \ar[ld]_p \ar[rd]^q & & \\ \mathbb{P} V_\Sigma\supset & \Sigma
 & & \Sigma^\vee & \subset\mathbb{P} V_\Sigma^\vee}
}

\medskip\noindent where we denoted by $I\subset \mathbb{P} V_\Sigma\times \mathbb{P} V_\Sigma^\vee$  
the incidence variety, defined as the set of pairs $(x,h)$ with $x\in\Sigma$ and $h$
a hyperplane containing the affine tangent space to $\Sigma$ at $h$. The variety 
$\Sigma^\vee$ parametrizing tangent hyperplanes to $\Sigma$ is its projective
dual variety. For all the cases we are interested in, the dual variety is non-degenerate,
which means that it is a hypersurface whose degree has been computed in the framework
of prehomogeneous vector spaces. 

Indeed, the variety $\Sigma$ is the closed orbit 
of the simple group $Aut(\Sigma)$ in a (projectivized) prehomogeneous space 
$\PP (V_\Sigma)$ whose orbit structure is known explicitly. Part of 
the structure is common to all cases: there is of course an open orbit, 
its complement is an irreducible hypersurface, and this hypersurface
contains an orbit closure $W_\Sigma$ whose complement is again a single
orbit. This can be observed case by case. 
 
\smallskip
If $\Sigma$ is the Cayley plane,  there are only three orbits in 
$\PP (V_\Sigma)$, the variety $\Sigma$ itself, its complement
in its secant variety, which is a cubic hypersurface, and the 
complement of this hypersurface. In particular $W_\Sigma=\Sigma$ in that case. 
An equation of the $E_6$-invariant 
cubic was first written down by Elie Cartan (in terms of the geometry 
of the $27$ lines on a smooth cubic  surface) and we therefore call 
it the Cartan cubic.   
 
If $\Sigma=\SS_{12}$, there are four orbits \cite{igusa},  whose closures are the 
whole projective space, a quartic hypersurface, $\Sigma$ itself, and 
an intermediate codimension seven variety which is our $W_\Sigma$. 

If $\Sigma=G(2,10)$, embedded in the projectivization of the space
of skew-symmetric forms, there are five orbits determined by the ranks of the 
forms; the invariant hypersurface is defined by the Pfaffian; 
$W_\Sigma$ is the space of forms of rank six at 
most, its  codimension is six. 

Finally, if $\Sigma=\SS_{14}$, there are exactly nine orbits which were 
first described by Popov \cite{popov}.
His results show that $W_\Sigma$ has codimension five. 

\smallskip
In each case we expect  $W_\Sigma$ to coincide with the singular locus of
the invariant hypersurface. At least we know that it contains this 
singular locus, just because of the orbit structure. 
 
The orbit structure is the same in the dual projective space $\PP (V_\Sigma^\vee)$,
in particular the invariant hypersurface is the projective dual variety of 
$\Sigma$. We summarize the relevant data in the following statement. 
We have no explanation, even conjectural, of these curious numerical 
coincidences.

\begin{lemm}\label{numbers}
The projective dual $\Sigma^\vee$ to $\Sigma$ is a hypersurface of 
degree $c-1$, whose singular locus has codimension at least $r-c+2$.
\end{lemm}

The following table gives $V_\Sigma$ for each case 
and the dual hypersurface $\Sigma^\vee$.

$$\begin{array}{lll}
\Sigma & V_\Sigma & \Sigma^\vee\\
\OP2& H_3(\mathbb{O}) & \mathrm{Cartan\; cubic} \\
\SS_{12} & \Delta_+  & \mathrm{Igusa\; quartic} \\
G(2,10) & \wedge^2\CC^{10} & \mathrm{Pfaffian\; quintic} \\
\SS_{14} & \Delta_+ & \mathrm{Popov\; octic}
\end{array}$$

\medskip 

Since the quotient 
of $V_\Sigma$ by the affine tangent bundle is the twisted normal bundle $N(-1)$, 
we have a natural identification $I\simeq \mathbb{P}(N(-1)^\vee)$. Moreover this 
identification is such that the relative tautological bundle $\mathcal{O}_{N(-1)}(1)$
coincides  with $q^*\mathcal{O}(1)$. 

Now, the duality theorem (see e.g. \cite{gkz}) asserts that $(\Sigma^\vee)^\vee =
\Sigma$, and moreover, that the diagram above is symmetric in the sense that $I$ can be 
defined as parametrizing the tangent hyperplanes to the hypersurface $\Sigma^\vee$ 
(at least at  smooth points, and then one needs to take the Zariski closure). 
If we denote by $F_\Sigma\in Sym^{c-1}V_\Sigma$ an equation of $\Sigma^\vee$, this 
means that $I$ can be obtained as the closure of the set of points $([dF_\Sigma(\xi)],[\xi])$,
where $\xi$ belongs to the cone over  $\Sigma^\vee_{reg}$. This means in particular
that 
we have a non-zero section 
$$\delta F_\Sigma\in H^0(I, Hom(q^*\mathcal{O}(-c+2),p^*\mathcal{O}(-1))),$$
vanishing precisely over the singular locus of $\Sigma^\vee$. Pushing this
section forward to $\Sigma$, we get a section 
$$\bar{\delta} F_\Sigma\in H^0(\Sigma, Sym^{c-2}(N(-1))(-1))=H^0(\Sigma, Sym^{c-2}N(-c+1))).$$
Now we can use the normal exact sequence $$0\ra T\Sigma\ra T\mathbb{P}(V_\Sigma)\otimes 
\mathcal{O}_\Sigma\ra N\ra 0.$$ Taking a suitable wedge power and twist, we get 
a long exact sequence 
$$0\ra \wedge^{c-2}T\Sigma (-c+1)\ra\cdots \ra Sym^{c-2}N(-c+1))\ra 0.$$
A careful check would then allow to conclude that this long exact sequence 
induces an isomorphism
$$H^0(\Sigma, Sym^{c-2}N(-c+1)))\simeq H^{c-2}(\Sigma, \wedge^{c-2}T\Sigma (-c+1)),$$
and the latter is Serre dual to $H^{r+2}(\Sigma,\Omega_\Sigma^{c-2}(c-r-1))$. 
In particular, we see that that our non-zero section $\bar{\delta} F_\Sigma$ induces a non-zero
linear form on this one dimensional cohomology group.

\subsection{Linear sections}

Let $X$ be a general linear section of $\Sigma$, of codimension
$s=r-c+1$. Then $X$ has dimension $n=2c-1$ and index $c-1$. 
The following statement is a special case of a result of Borcea \cite{borcea}:

\begin{lemm}
Any small deformation of $X$ is a linear section of $\Sigma$.
\end{lemm}

In particular, the number of moduli for $X$ is the dimension $m$ 
of the quotient of the Grassmannian $G(s,V_{\Sigma})$ by the 
simple group $Aut(\Sigma)$, of dimension $\delta$. The relevant 
data appear in the table below. 

$$\begin{array}{cccccc}
\Sigma  & s & N & Aut(\Sigma) & \delta & m\\
\OP2 & 9 &  26 & E_6 & 78 & 84 \\
\SS_{12} & 6 &  31 & Spin_{12} & 66 & 90 \\
G(2,10) & 5 & 44 & PSL_{10} & 99 & 101\\
\SS_{14} & 4 & 63 & Spin_{14} & 91 & 149
\end{array}$$

\begin{theo}\label{h=m}
A generic linear section $X$ of $\Sigma$ has primitive 
cohomology only in middle degree. 
Its non-zero Hodge numbers $h^{p,n-p}$, 
for $p\ge c$, are the following:
$$h^{c+1,c-2}=1, \qquad h^{c,c-1}=m.$$
\end{theo}

\proof Let us denote by $L\subset V_{\Sigma}^\vee$ the space 
of linear forms that defines $X$. Its dimension is $s=r-c+1$.
We use the conormal sequence
$$0\ra \cO_X(-1)\otimes L\ra\Omega_{\Sigma |X}^1\ra\Omega_X^1\ra 0,$$
and its wedge powers, for  positive integers $p$, 
$$0\ra \cO_X(-p)\otimes Sym^pL
\ra\cdots\ra
\Omega_{\Sigma |X}^{p-1}(-1)\otimes L\ra\Omega_{\Sigma |X}^{p}
\ra \Omega_X^p\ra 0.$$
We claim that for $p\le c-1$, this induces an isomorphism 
$$H^q(X,\Omega_{\Sigma |X}^{p})\simeq H^q(X,\Omega_X^p).$$
To check this, we will prove that the other vector bundles involved 
in the long exact sequence above, are all acyclic. These bundles are
the bundles of twisted forms $\Omega_{\Sigma |X}^{p-k}(-k)$, for $1\le k\le p$. 
For example, for $k=p$ we observe that 
$$H^q(X,\cO_X(-p))=0 \qquad \forall q\le n-1$$
by Kodaira's vanishing theorem, while by Serre duality 
$$H^q(X,\cO_X(-p))=H^{n-q}(X,\cO_X(p-c+1))^\vee=0\qquad \forall q\ge 1$$
when $p\le c-2$. For the remaining terms, we use the Koszul resolution 
of the structure sheaf of $X$, 
$$0\ra\cO_\Sigma(c-r-1)\otimes \wedge^{r-c+1}L
\ra\cdots\ra \cO_\Sigma(-1)\otimes L\ra\cO_\Sigma\ra\cO_X\ra 0.$$
Twisting this long exact sequence by $\Omega_{\Sigma}^{p-k}(-k)$,
we see that in order to prove the acyclicity of $\Omega_{\Sigma |X}^{p-k}(-k)$,
it is enough to check the acyclicity of $\Omega_{\Sigma}^{p-k}(-k-\ell)$
for $0\le\ell\le r-c+1$. But this follows from Lemma \ref{van}.

By the same Lemma, 
 $\Omega_{\Sigma}^{p}(-\ell)$ is acyclic
for $1\le\ell\le r-c+1$ and $p\le c-3$. This implies that 
$$H^q(\Sigma,\Omega_{\Sigma}^{p})\simeq H^q(X,\Omega_{\Sigma |X}^{p})
\simeq H^q(X,\Omega_X^p)$$
for all $q$ and for all $p\le c-2$. 

For $p=c-2$ there is a unique non acyclic bundle involved, 
$\Omega_{\Sigma}^{c-2}(c-r-1)$: by Lemma \ref{nonvan}, this 
bundle has a one-dimensional cohomology group in degree $r+2$, 
but the other cohomology groups vanish. Then the 
Koszul resolution of $\cO_X$ implies that 
$$H^q(X,\Omega_{\Sigma |X}^{c-2})_0=
H^{q+r-c+1}(\Sigma,\Omega_{\Sigma}^{c-2}(c-r-1))=\delta_{q,c+1}\CC.$$

Finally, for $p=c-1$ we find two non-zero cohomology groups
for the bundles involved in the complex
$$0\ra \cO_X(1-c)\otimes Sym^{c-1}L
\ra\cdots\ra
\Omega_{\Sigma |X}^{c-2}(-1)\otimes L\ra\Omega_{\Sigma |X}^{c-1}
\ra \Omega_X^{c-1}\ra 0.$$
Indeed, by Lemma \ref{van} $\Omega_{\Sigma}^{c-2}(c-r)$ is acyclic. 
Then, by the same argument as above, the cohomology of the bundle 
$\Omega_{\Sigma |X}^{c-2}(-1)$ is 
$$H^q(X,\Omega_{\Sigma |X}^{c-2}(-1))=
H^{q+r-c}(\Sigma,\Omega_{\Sigma}^{c-2}(c-r-1)\otimes \wedge^{r-c}L)=
\delta_{q,c+2}L^\vee.$$
On the other hand, the line bundle $\cO_X(1-c)$ is the canonical 
bundle of $X$, and therefore it has non trivial cohomology in degree
$n=2c-1$, and only in that degree. Note that $\Omega_{\Sigma |X}^{c-1}$
also has non-zero cohomology, but only in degree $c-1$, fo which we get 
the same cohomology group as for  $\Omega_{\Sigma}^{c-1}$. We deduce 
a long exact sequence 
\begin{eqnarray*}
 \cdots\ra H^q(X,\Omega_X^{c-1})_0\ra
H^{q+c-1}(X,\omega_X)\otimes Sym^{c-1}L
\ra \hspace*{3cm} {\it } \\
 \hspace*{2cm}\ra H^{q+2}(X,\Omega_{\Sigma |X}^{c-2}(-1))\otimes L\ra 
H^{q+1}(X,\Omega_X^{c-1})_0\ra\cdots
\end{eqnarray*}
But note that this long exact sequence has very few non-zero terms,
appearing precisely for $q=c$. Note also that by Hodge symmetry, we 
already know that $H^{c+1}(X,\Omega_X^{c-1})=H^{c}(X,\Omega_X^{c-2})=0$. 
What remains is the short exact sequence
\begin{equation}\label{moduli}
0\ra H^{c}(X,\Omega_X^{c-1})\ra Sym^{c-1}L\ra 
L^\vee\otimes L\ra 0,
\end{equation}
from which we can easily compute $h^{c,c-1}$. And we conclude, as claimed, 
that $h^{c,c-1}=m$. \qed

\medskip\noindent {\it Remark}. In general, if $\Sigma$ has a non degenerate dual of
degree $d\le n$, we could expect $H^{n+1-d}(\Sigma,\Omega_\Sigma^{d-1}(d-r))$ to be 
non-zero, and even one-dimensional. Under favorable circumstances, this
should imply that if $X$ is a general linear section of codimension $r-d$, hence
dimension $c+d$, then $h^{c+1}(\Omega_X^{d-1})=1$.

\subsection{Deformations}
Observe that the conclusion of Theorem \ref{h=m} has been obtained 
by a case by case computation, as a coincidence between two numbers
that were computed from quite different, and at first sight, unrelated
 data. 

On one hand, the number  $m$ of moduli has been computed as the difference between the dimension 
of a suitable Grassmannian and that of the automorphism group of $\Sigma$. 
Indeed, $H^1(X,TX)$ is given by the exact sequence 
$$0\ra H^0(X,T\Sigma_{|X})\ra H^0(X,L^\vee(1))\ra H^1(X,TX)\ra 0.$$
It is straightforward to check that 
$$H^0(X,T\Sigma_{|X})=H^0(\Sigma,T\Sigma)=aut(\Sigma),$$
the Lie algebra of $Aut(\Sigma)$, 
while $H^0(X,L^\vee(1))=L^\vee\otimes V_\Sigma^\vee/L$ is the tangent 
space to the Grassmannian of subspaces of $V_\Sigma$, at its point
defined by $L$. This tangent space is also the quotient of $gl(V_\Sigma)$
by the stabilizer of $L$, hence the identification 
$$H^1(X,TX)\simeq gl(V_\Sigma)/(aut(\Sigma)+stab(L)).$$
On the other hand, the exact sequence (\ref{moduli}) shows that 
$$H^c(X,\Omega^{c-1}_X)^\vee\simeq Sym^{c-1}L^\vee/gl(L).$$ 
As we have seen, the connection between these two descriptions is
provided by the equation $F_\Sigma$, of degree $c-1$, of  the 
dual hypersurface $\Sigma^\vee$. We have a natural map from
$gl(V_\Sigma)$ to $Sym^{c-1}V_\Sigma$, sending $u$ to $F_\Sigma
\circ u$, then to $Sym^{c-1}L^\vee$ by restriction to $L$. 
Since $aut(\Sigma)$ kills $F_\Sigma$, this induces a map 
\begin{equation}\label{diff}
gl(V_\Sigma)/(aut(\Sigma)+stab(L))\ra 
 Sym^{c-1}L^\vee/gl(L),
\end{equation}
which is an isomorphism if and only if $X$ is of Calabi-Yau
type. 

\medskip
This has a modular interpretation: we can associate to $X$ the 
hypersurface $X^\vee = \PP L\cap \Sigma^\vee$ (which is of course {\it not}
the projective dual variety of $X$). By Lemma \ref{numbers}, this hypersurface
will be smooth of degree $c-1$ for a generic $L$. Moreover, if $X^\vee$ has discrete automorphism group, 
the quotient $Sym^{c-1}L^\vee/gl(L)$ can be interpreted as $H^1(X^\vee,TX^\vee)$,
Otherwise said, the correspondence between $X$ and $X^\vee$ defines a map
from the moduli space of linear  sections of $\Sigma$ to the moduli
space of degree $c-1$ hypersurfaces in $\PP^{r-c}$, and the morphism (\ref{diff})
is nothing else than the differential of that map. We deduce:

\begin{prop}
The variety $X$ is a Fano manifold of Calabi-Yau type if and only 
if its small deformations induce a versal family of degree $c-1$ 
hypersurfaces in $\PP^{r-c}$. 
\end{prop}

Note that if this is true, we can conclude that a generic hypersurface of degree $c-1$ 
in $\PP^{r-c}$ can be defined, up to isomorphism, as a linear section of $\Sigma^\vee$ 
in a finite (non-zero) number of different ways. 

For the quintic threefold, it has already been noticed by Beauville that there
exists finitely many Pfaffian representations (see Schreyer's Appendix to \cite{beau}). 
Beauville asked how many such representations do exist : in principle the answer 
is given by some Donaldson-Thomas invariant. We can ask the same question
in our three other cases. The case of cubic sevenfolds is treated in \cite{im-cub}
where we prove that the very same phenomena hold. That's also what we expect 
for the two remaining cases, which will deserve further investigations.

\medskip
It is also quite remarkable that we can associate to $X$ either a Calabi-Yau 
threefold,
or another Fano manifold of Calabi-Yau type,   
by considering either the dual $X^\vee$, or a double covering of $\PP^{r-c}$,
branched over $X^\vee$. We get the following varieties $X^*$:

$$\begin{array}{ll}
X & X^* \\
\OP2\cap\PP^{17} & \mathrm{cubic\; sevenfold} \\
\SS_{12}\cap\PP^{25} &   \mathrm{double\; quartic\; fivefold} \\
G(2,10)\cap\PP^{40} & \mathrm{quintic\; threefold} \\
\SS_{14}\cap\PP^{59} & \mathrm{double\; octic\; threefold}
\end{array}$$

\medskip
The last two varieties are famous examples of Calabi-Yau threefolds. 
The first two were among our very first examples of Fano manifolds 
of Calabi-Yau type. They are not Calabi-Yau's stricto sensu but, 
considered as (compactified) Landau-Ginzburg models, we have seen that 
they appeared in the litterature
as mirrors to certain rigid Calabi-Yau threefolds \cite{cdp,sch}. 

\smallskip
Note that if the resulting map $\mathcal{K}_X\ra \mathcal{K}_{X^*}$ from
the base of the Kuranishi family of $X$, to that of $X^*$, is etale, 
then we have the same property for the gauged families 
and we can pull-back the integrable system in intermediate Jacobians 
of $X^*$ and its deformations (see Theorem \ref{int}). Although we did not 
check this, it is quite likely that we should recover the integrable 
system in intermediate Jacobians of $X$ itself and its deformations. 
In particular the intermediate Jacobians of $X$ and $X^*$ should
be isomorphic. 

\subsection{Homological projective duality}

We expect that our series of examples of Fano manifolds of 
Calabi-Yau type should give rise to interesting instances 
of homological projective duality \cite{kuz2}. More precisely, 
the derived categories of $X$ and $X^*$ should both contain
a three-dimensional Calabi-Yau category, say $\mathcal{A}_X\subset
\mathcal{D}^b(Coh(X))$ and  $\mathcal{A}_{X^*}\subset
\mathcal{D}^b(Coh(X^*))$ and there should exist a natural
equivalence between $\mathcal{A}_X$ and $\mathcal{A}_{X^*}$.

For $X^*$ this was observed by Kuznetsov:

\begin{prop}   
Let $Y$ be among a the following types of Fano manifolds of Calabi-Yau type:
\begin{enumerate}
\item a cubic sevenfold,
\item a cubic hypersurface in a six-dimensional quadric,
\item a double-cover of $\PP^5$ branched over a smmoth quartic. 
\end{enumerate}
Then the derived category $\cD^b(Y)$ contains a natural subcategory 
$\cA_Y$ which is a three-dimensional Calabi-Yau category. 
\end{prop}

\proof For any Fano hypersurface $Y\subset\PP^{n+1}$ of degree $d$, and
index $\iota_Y=n+2-d$, the collection 
$$\langle \cO_Y,\ldots ,\cO_Y(\iota_Y-1)\rangle$$
is obviously exceptional. Kuznetsov proved in \cite{kuz1}, Corollary 4.3,
that if $d$ divides $n+2$, its left orthogonal $\cA_Y$ is a Calabi-Yau category of
dimension $n-2\iota_Y/d$. In particular, for $d=3$ and $n=7$, we can conclude that
$\cA_Y$ is a three-dimensional  Calabi-Yau category. 

A very similar statement holds for double covers. Suppose that $Y$ is a double-cover
of $\PP^n$ branched over a general hypersurface of degree $2d$. Then $Y$ is Fano 
for $d\le n$, and $\omega_Y^{-1}=f^*\cO_{\PP^n}(n+1-d)$, so that $\iota_Y=n+1-d$. 
Here, for the same reason as before,  the collection $$\langle \cO_Y,f^*\cO_{\PP^n}(1),\ldots ,
f^*\cO_{\PP^n}(\iota_Y-1)\rangle$$ is exceptional. Kuznetsov proved that its
orthogonal is, when $d$ divides $n+1$, a Calabi-Yau category of dimension 
$n-\iota_Y/d$. For $d=2$ and $n=5$, we can conclude that
$\cA_Y$ is a three-dimensional  Calabi-Yau category. 

Finally the case of hypersurfaces in quadrics in very similar. In dimension $6$,
the derived category has a Lefschetz decomposition  $\cD^b(\QQ^6) = \langle A , A(3) \rangle,$
where 
$$A = \langle  S , O_{\QQ^6} , O_{\QQ^6}(1), O_{\QQ^6}(2)\rangle ,$$ 
and $S$ denotes the rank four spin bundle. If $i : Y\hookrightarrow \QQ^6$ is a cubic hypersurface, 
the functor $i^* : D^b(\QQ^6) \ra D^b(Y)$ is fully faithful on $A$
and the orthogonal $\cA_Y$ to $i^*(A)$ in $D^b(Y)$ is a three-dimensional
Calabi-Yau category.\qed

\smallskip
The appearance of these three dimensional {\it non-commutative} 
Calabi-Yau's in connection with our manifolds of Calabi-Yau type
(but not three dimensional!) is a quite remarkable phenomenon,
which would deserve to be better understood.

\section{Appendix: Proof of Lemma \ref{nonvan}}

First observe that, $\Sigma$ being a Hermitian symmetric space, its 
cotangent bundle is an  irreducible homogeneous bundle. In particular
its cohomology groups, and that of its twists as well, can be computed using 
Bott's theorem. Moreover the same is true for any bundle of $p$-forms
on $\Sigma$, since it must split into a direct sum of irreducible homogeneous 
bundles. 

Recall that an irreducible homogeneous vector bundle on a homogeneous  
space $\Sigma=G/P$ is defined by an irreducible representation of $P$, 
which is in turn defined by its highest weight: some weight $\lambda$ 
of $G$ which is $P$-dominant. We denote the corresponding bundle by 
$\cE_\lambda$.

When $\Sigma$ is a Hermitian symmetric space with Picard number one, 
$P$ is a maximal parabolic subgroup defined by some simple root
$\alpha$, and the cotangent bundle is simply 
$$\Omega^1_\Sigma=\cE_{-\alpha}.$$
More generally, a combinatorial formula has been obtained by Kostant 
for the irreducible components of each $\Omega^p_\Sigma$ \cite{kos}.

Bott's theorem for an irreducible  homogeneous vector bundle $\cE_{\lambda}$
can be stated as  the following recipe. 
Add $\rho$, the half-sum of the positive roots of $G$, to $\lambda$. If 
$\langle\lambda+\rho,\alpha^\vee\rangle=0$ for some coroot $\alpha^\vee$, then
$\cE_{\lambda}$ is acyclic. Otherwise, there is a unique $w$ in the Weyl 
group $W$ of $G$ such that $w(\lambda+\rho)=\mu+\rho$ for some dominant
weight $\mu$. Then $\cE_{\lambda}$ has a unique non-zero cohomology 
group, in degree $\ell(w)$, and it is an irreducible $G$-module of 
lowest weight $-\mu$. In particular, it is one dimensional if and only 
if $\mu=0$. 

\medskip\noindent 
{\it First case : $\Sigma=G(2,10)$}. Here $\Omega^1_\Sigma=Q^*\otimes E$,
where $E$ denotes the tautological rank two bundle and $Q$ the quotient 
rank eight bundle. The claim is that 
 $$H^{12}(\Sigma,\Omega_\Sigma^{4}(-5))=\CC.$$
The decomposition of $\Omega_\Sigma^{4}$ into irreducible components is 
$$\Omega_\Sigma^{4}=\wedge^4Q^*\otimes S^4E\oplus S_{211}Q^*\otimes S_{31}E
\oplus S_{22}Q^*\otimes S_{22}E.$$
We apply Bott's theorem to each factor. The weights of the three factors
are 
$$\begin{array}{l}
(0,0,0,0,-1,-1,-1,-1,4,0),\\ 
(0,0,0,0,0,-1,-1,-2,3,1),\\ 
(0,0,0,0,0,0,-2,-2,2,2).
\end{array}$$ 
Twisting by $\cO_\Sigma(-5)$ and adding $\rho$ we get 
$$\begin{array}{l} 
(4,3,2,1,-1,-2,-3,-4,5,0), \\
(4,3,2,1,0,-2,-3,-5,4,1), \\ 
(4,3,2,1,0,-1,-3,-4,3,2).
\end{array}$$
The last two sequences have repeated entries, indicating 
that the corresponding bundles are acyclic. But the first one has pairwise
distinct entries, forming, once ordered, a consecutive sequence. This indicates
that the corresponding bundle has one non-zero cohomology group, of 
dimension one, appearing in degree equal to the number of inversions of the
sequence, which is twelve. This proves the claim.  

\medskip\noindent 
{\it Second case : $\Sigma=\SS_{12}$}. Here $\Omega^1_\Sigma=\wedge^2 E$,
where $E$ denotes the tautological rank six bundle. The claim is that
$$H^{12}(\Sigma,\Omega_\Sigma^{3}(-6))=\CC.$$
The decomposition of $\Omega_\Sigma^{3}$ into irreducible components is
$$\Omega_\Sigma^{3}=S_{411100}E\oplus S_{222000}E.$$
Twisting by $\cO_\Sigma(-6)$ and adding $\rho=(5,4,3,2,1,0)$ we get 
the weights $(2,1,0,-3,-4,-5)$ and $(2,0,-1,-2,-3,-7)$. The last 
sequence has two opposite entries, indicating that the corresponding bundle
is acyclic. But this is not the case of the first sequence, whose 
absolute values form a consecutive sequence starting from zero. This 
implies that the corresponding bundle has a unique non-zero cohomology group,
which is one dimensional. It appears in degree equal to the number of 
pairs of entries of the sequence whose sum is negative; this is twelve 
and the claim follows. 

  \medskip\noindent 
{\it Third case : $\Sigma=\SS_{14}$}. Here again $\Omega^1_\Sigma=\wedge^2 E$,
where $E$ denotes the tautological rank seven bundle. The claim is that
$$H^{14}(\Sigma,\Omega_\Sigma^{7}(-4))=\CC.$$
The decomposition of $\Omega_\Sigma^{7}$ into irreducible components is
$$\Omega_\Sigma^{7}=S_{6221111}E\oplus S_{5322110}E
\oplus S_{4333100}E\oplus S_{4422200}E.$$
Twisting by $\cO_\Sigma(-7)$ and adding $\rho=(6,5,4,3,2,1,0)$ we get 
the weights
$$\begin{array}{l} 
(3,2,1,0,-2,-3,-8), \\
(4,3,1,-2,-3,-4,-6), \\ 
(4,2,1,-1,-2,-4,-7)), \\
(4,3,0,-1,-2,-5,-6).
\end{array}$$
The first three have pairs of opposite entries, which implies that
the corresponding bundles are acyclic. This is not the case of the
last sequence, whose absolute values form a consecutive sequence 
starting from zero. As in the previous case, this implies the claim.

  \medskip\noindent 
{\it Fourth case : $\Sigma=\OP2$}. For convenience we encode
each irreducible vector bundle $\cE_{\lambda}$ on the Cayley plane 
$\OP2$ by a weighted Dynkin diagram of type $E_6$, whose weights are the 
coordinates of $\lambda$ on the basis of fundamental weights (themselves 
in natural bijection with the nodes of the Dynkin diagram). For example
we have 
$$
\Omega^1_\Sigma \simeq \stackrel{-21000}{\;\;\scriptstyle 0},\qquad
\Omega^2_\Sigma  \simeq \stackrel{-30100}{\;\;\scriptstyle 0}.$$
 
The claim to be checked is that 
$$H^{14}(\Sigma,\Omega_\Sigma^{2}(-9))=\CC.$$
The root system, and the Weyl group $W$ of $E_6$ being rather complicated, 
we will use a different strategy than in the other case. 
If $\lambda=-3\omega_1+\omega_4$ denotes the highest 
weight of $\Omega^2_\Sigma$, what we need
to check is that $\mu=\lambda-9\omega_1+\rho$ belongs to the $W$-orbit
of $\rho$. Recall that $\rho$, the half sum of the positive roots,
is also the sum of the fundamental weights. We can thus write
$$\rho = \stackrel{11111}{\scriptstyle 1},\qquad \mu = \stackrel{-B 1211}{\;\;\scriptstyle 1},$$
where $B$ stands for eleven (and $A$ will stand for ten, below). 
Now we can let $W$ act through the simple reflections $s_i$ associated
to the simple root $\alpha_i$. The recipe for the action of $s_i$ on a 
weight $\omega$, written in the basis of fundamental weights, is as follows:
add the $i$-th coordinate to the $j$-th if $j$ is connected to $i$ in the
Dynkin diagram; then change the $i$th coordinate into its opposite. 
Here we start with a weight having a negative coordinate. Applying 
successively simple reflections corresponding to nodes with negative
coordinates, we will end up with a dominant weight. Starting from
$\mu$, the game goes as follows. 
$$\begin{array}{l}
\stackrel{-B 1211}{\;\;\scriptstyle 1} \quad\ra\quad \stackrel{B -A211}{\;\;\;\;\scriptstyle 1} 
\quad\ra\quad \stackrel{1A-811}{\;\;\scriptstyle 1}
\quad\ra\quad \stackrel{128-71}{\scriptstyle -7}\quad\ra\quad \stackrel{121-71}{\scriptstyle 7}\quad\ra \\
  \\
\quad\ra\quad \stackrel{12-67-6}{\scriptstyle 7}\quad\ra\quad \stackrel{12-616}{\scriptstyle 7}
\quad\ra\quad \stackrel{1-46-56}{\;\scriptstyle 1}\quad\ra\quad \stackrel{1-4151}{\;\;\scriptstyle 1} 
\quad\ra\quad \stackrel{-34-351}{\;\;\scriptstyle 1}\quad\ra \\
   \\
\quad\ra\quad \stackrel{31-351}{\;\;\scriptstyle 1}
\quad\ra\quad \stackrel{3-2321}{\scriptstyle -2}\quad\ra\quad \stackrel{3-2121}{\;\;\scriptstyle 2} 
\quad\ra\quad \stackrel{12-121}{\;\;\scriptstyle 2}\quad\ra\quad \stackrel{11111}{\;\scriptstyle 1}
\end{array}$$

\smallskip
So we end up with $\rho$ after $14$ operations. This indicates that 
we get a one dimensional cohomology group in degree fourteen. The claim
is proved. \qed

\medskip 
With the very same method we can prove that the one dimensional
cohomology group whose existence is asserted by Lemma \ref{nonvan},
is the only non-zero cohomology group of twisted form in a wide range. 
The precise statement is the following:

\begin{lemm}\label{van}
Suppose that $p\le c-1$ and $1\le k\le r-p$. Then 
$$H^{q}(\Sigma,\Omega_\Sigma^{p}(-k))=\delta_{q,r+2}
\delta_{p,c-2}\delta_{k,r-c+1}\CC.$$
\end{lemm}

\end{document}